\documentclass[reqno]{amsart}

\newtheorem{theorem}{Theorem}

\theoremstyle{definition}

\usepackage{graphicx}
\usepackage{pstricks, enumerate, pst-node, pst-text, pst-plot}

\theoremstyle{remark}
\newtheorem{remark}[theorem]{Remark}

\numberwithin{equation}{section}

\newcommand{\intav}[1]{\mathchoice {\mathop{\vrule width 6pt height 3 pt depth  -2.5pt
\kern -8pt \intop}\nolimits_{\kern -6pt#1}} {\mathop{\vrule width
5pt height 3  pt depth -2.6pt \kern -6pt \intop}\nolimits_{#1}}
{\mathop{\vrule width 5pt height 3 pt depth -2.6pt \kern -6pt
\intop}\nolimits_{#1}} {\mathop{\vrule width 5pt height 3 pt depth
-2.6pt \kern -6pt \intop}\nolimits_{#1}}}

\newcommand{\intavl}[1]{\mathchoice {\mathop{\vrule width 6pt height 3 pt depth  -2.5pt
\kern -8pt \intop}\limits_{\kern -6pt#1}} {\mathop{\vrule width 5pt
height 3  pt depth -2.6pt \kern -6pt \intop}\nolimits_{#1}}
{\mathop{\vrule width 5pt height 3 pt depth -2.6pt \kern -6pt
\intop}\nolimits_{#1}} {\mathop{\vrule width 5pt height 3 pt depth
-2.6pt \kern -6pt \intop}\nolimits_{#1}}}

 \newcommand{\mc}{\mathcal}

 \newcommand{\F}{\mc{F}}

 \newcommand{\U}{\mc{U}}

 \newcommand{\C}{\mathbb{C}}
 \newcommand{\R}{\mathbb{R}}

 \newcommand{\Sp}{\mathbb{S}}

 \newcommand{\dx}{\text{\rm d}x}

 \newcommand{\Kahler}{K\"{a}hler}
 \newcommand{\Rc}{\text{Rc}}
 \newcommand{\ab}{\alpha\bar{\beta}}
 \newcommand{\Czero}{\mathbb{C}^{2}\backslash\{0\}}
 \newcommand{\Cnzero}{\mathbb{C}^{n}\backslash\{0\}}

 \newcommand{\CP}{\mathbb{CP}}
 
 \newcommand{\al}{\alpha}
 \newcommand{\be}{\beta}
 \newcommand{\bebar}{\bar{\beta}}
 
 \newcommand{\dalphabetabar}{\dfrac{\partial^2 }{\partial z^{\alpha} \partial \bar{z}^{\beta}}}
 \usepackage{epsfig}

\begin{document}

\title[Non-negative Ricci under Ricci flow]{Non-negative Ricci curvature on closed manifolds under Ricci flow}

\author[Davi M\'aximo]{Davi M\'aximo}
\address{Department of Mathematics, University of Texas at Austin, TX 78712, USA.}
\email{maximo@math.utexas.edu}

\


\date{October, 14, 2009}


\begin{abstract}
In this short note we show that non-negative Ricci curvature is not preserved under Ricci flow for closed manifolds of dimensions four and above, strengthening a previous result of Knopf in \cite{K} for complete non-compact manifolds of bounded curvature. This brings down to four dimensions a similar result B\"{o}hm and Wilking have for dimensions twelve and above, \cite{BW}. Moreover, the manifolds constructed here are \Kahler\, manifolds and relate to a question raised by Xiuxiong Chen in \cite{XC}, \cite{XCL}. 
\end{abstract}
\maketitle

\section{Introduction}
\indent Ricci flow is a flow of Riemannian metrics designed to improve a given initial Riemannian metric and has been a powerful tool in understanding geometry and topology of manifolds. A nice and important feature of this flow is that it preserves the ``non-negativity" of various curvature conditions. \\
\indent The flow was introduced on closed manifolds by Richard Hamilton in \cite{H1}. He proved that non-negative scalar curvature is preserved by the flow in all dimensions. He also showed that the flow preserves non-negative Ricci and non-negative sectional curvatures on three dimensional closed manifolds. Later, in \cite{H2}, Hamilton proved that Ricci flow also preserves the non-negativity of the Riemann curvature tensor in all dimensions. Also in higher dimensions, independent work of Huisken, Mar\-ge\-rin and Nishikawa (see \cite[Section 7.1]{CNL}) demonstrated the invariance of certain sets defined by pinching conditions. Haiwen Chen, in \cite{HC}, showed that 2-nonnegativity of the Riemann tensor is preserved in all dimensions. In dimension four, Hamilton showed in \cite{H3} that positive isotropic curvature  is invariant under Ricci flow and it has recently been proved independently by Huy T. Nguyen \cite{Ngu} and by Brendle and Schoen \cite{BS} that this holds in dimensions five and above. In their important paper \cite{BS}, Brendle and Schoen also proved that each of the properties of $M^n\times \R$ and $M^n\times\R^2$ having non-negative isotropic curvature is preserved by Ricci flow for all $n\geq 4$; critically, the latter property is weaker than pointwise $\frac{1}{4}$-pinched sectional curvature, but stronger than non-negative sectional curvature. In related work, Brendle found that if $M^n \times \Sp^2(1)$ has nonnegative isotropic curvature then this remains so when the metric on M is evolved by the Ricci flow (where $\Sp^2(1)$ denotes a two-sphere of constant curvature 1); we refer to the reader to \cite{Bre}. B\"{o}hm and Wilking in \cite{BW1} provided a method for constructing new invariant curvature conditions from known ones. Finally, during summer 2009, Burkhard Wilking announced a interesting result that provides a unified approach for understanding many previously discovered invariant curvature conditions. \\  
\indent Ricci flow was extended to complete non-compact manifolds with bounded curvature by Wan-Xiong Shi in \cite{Sh2} and to \Kahler\, manifolds by Huai-Dong Cao in \cite{Cao3} and by Shigetoshi Bando in \cite{B}. Similar invariant curvature conditions were later found for these flows too. In \cite{Sh1}, Shi proved that non-negative Ricci curvature is preserved by the flow on complete manifolds of three dimensions. Bando and Mok, in \cite{Mok} and \cite{B}, proved that non-negative bisectional curvature (which is a stronger assumption than and implies non-negative Ricci curvature, see \cite{GK}) on closed \Kahler\, manifolds is also an invariant curvature condition and this later was extended by Shi in \cite{Sh3} to complete non-compact \Kahler\, manifolds with bounded curvature. \\
\indent There has also been work showing that some curvature conditions are not preserved in general. In \cite{Ni}, Lei Ni exhibited complete non-compact Riemannian manifolds with bounded non-negative sectional curvature of dimensions four and above where the Ricci flow does not preserve the non-negativity of the sectional curvature. Dan Knopf later showed in \cite{K} examples of complete non-compact \Kahler\, manifolds of bounded curvature and non-negative Ricci curvature whose \Kahler-Ricci evolutions immediately acquire Ricci curvature of mixed sign, for each complex dimension above one. Finally, in \cite{BW}, B\"{o}hm and Wilking showed an example of a closed twelve dimensional manifold where Ricci flow evolves an initial homogeneous metric of positive sectional curvature to metrics with mixed Ricci curvature. By taking products with spheres, they settled that positive Ricci is not invariant under the flow for dimensions twelve and above.\\
\indent Starting on \cite{XCT1} and later with \cite{XCT2}, Xiuxiong Chen and Gang Tian answered affirmatively the following important question on Ricci flow: on a compact \Kahler-Einstein manifold, does the \Kahler-Ricci flow converge to a \Kahler-Einstein metric if the initial metric has positive bisectional curvature? They remarked that what their argument really needed was for the Ricci curvature to be positive along the \Kahler-Ricci flow. Since the positivity of Ricci curvature being preserved under the Ricci flow was unknown, they used the fact that the positivity of the bisectional curvature is preserved and that it implies positive Ricci. Since their work, others tried to remove or weaken this strong assumption of positive bisectional curvature and we refer to the interested reader \cite{PS1}, \cite{PS2}, \cite{XC}, \cite{XCL}. 
In $1988$, Richard Hamilton and Huai-Dong Cao announced a proof that positive orthogonal bisectional curvature is preserved along the flow, and that this suffices to guarantee that non-negative Ricci is preserved. More recently, Xiuxiong Chen published a proof in \cite{XC} that non-negative Ricci curvature is preserved on a \Kahler\,
solution, assuming its bisectional curvature remains positive as long as it exists. Moreover, under the same assumption, he proved that if the initial Ricci curvature is not positive, then the lower bound of Ricci will increase with time and actually reach zero as time goes to infinity. Motivated by his nice result in \cite{XC}, Chen asked later in \cite{XCL} if some form of lower bound of the Ricci curvature would be preserved under Ricci flow at least for closed \Kahler\, manifolds.\\ 
\indent The purpose of this short note is to show that this is {\it not} the case in general. For each complex dimension greater than one, we construct a closed \Kahler\, manifold with non-negative Ricci curvature whose evolution under Ricci flow immediately acquire Ricci curvature of mixed sign. Let $M$ denote $\CP^2$ blown up at two points (in the next section, following Calabi in \cite{C}, we realize $M$ as the total space of a bundle $\CP^1\hookrightarrow M \rightarrow\CP^1$; see \cite[Section 2.7.2]{Chow} for more details). We prove the following:

\begin{theorem}\label{theorem}
There exist \Kahler\, metrics on $M$ with non-negative Ricci curvature that immediately acquire mixed Ricci curvature when evolved by Ricci flow.
\end{theorem}

\begin{remark}
By taking products of $M$ with spheres, we lower B\"{o}hm and Wilking previous result in \cite{BW} to real dimension four and above. Observe that these products will not be \Kahler\, manifolds anymore.   
\end{remark}

\indent Theorem \ref{theorem} can actually be generalized to all complex dimensions above one:

\begin{theorem}\label{theorem2}
For $n\geq 2$, let $M_n$ be $\CP^n$ blown up at two points. Then, there exist \Kahler\, metrics on $M_n$ with non-negative Ricci curvature that immediately acquire mixed Ricci curvature when evolved by Ricci flow.
\end{theorem} 

\begin{remark}
Examples found in theorems \ref{theorem} and \ref{theorem2} can be thought as compactifications of the ones found by Knopf in \cite{K}.
\end{remark}

\begin{remark}
A three-manifold $M^3$ has non-negative Ricci curvature if and only if $M^3\times \R$ has non-negative isotropic curvature. Our result shows that non-negative Ricci does not behave well in higher dimensions, so Brendle and Schoen's observation that non-negative isotropic curvature on $M^n\times\R$ is invariant suggests that this property is the ``right'' generalization to higher dimensions.
\end{remark}

\indent The rest of this paper is organized as follows. We first deal with the critical case $n=2$ and begin by proving Theorem \ref{theorem} in Sections $2$ and $3$: in $2$, we review $U(2)$-invariant metrics on $\Czero$ and construct a initial metric on $M$ with non-negative Ricci curvature, and in Section $3$, the proof of Theorem \ref{theorem} is finished. Finally, in Section $4$, we discuss the general case stated on Theorem \ref{theorem2}.   

\section{An Initial Metric of Non-Negative Ricci Curvature} 
In this section we consider rotationally symmetric \Kahler\, metrics on $\Czero$ and derive an initial metric of non-negative Ricci curvature on the twisted line bundle 
$\CP^1\hookrightarrow M \rightarrow\CP^1$ of Calabi \cite[Section 3]{C}. Metrics of this sort were also used by Cao in \cite{Cao1, Cao2} and Feldman, Ilmanen $\&$ Knopf in 
\cite{FIK} to construct examples of gradient Ricci solitons. For convenience of the reader, we start following Calabi in \cite{C} and realizing $M_n$ as the total space 
of a bundle $\CP^1\hookrightarrow M \rightarrow\CP^{n-1}$ (for more details about this construction, see \cite[Section 2.7.2]{Chow}). 

\subsection{Calabi's Bundle Construction} Cover the projective space $\CP^{n-1}$ with its usual $n$ charts $\varphi_{\alpha}:\U_\alpha\longrightarrow \C^{n-1}$, where 
$\U_{\alpha}=\{[x_1:x_2:\cdots:x_n]\in\CP^{n-1}:x_\alpha\neq0\}$ and $\varphi_\alpha:[x_1:x_2:\cdots:x_n]\mapsto\left(\frac{x_1}{x_\alpha},\ldots,
\frac{x_{\alpha-1}}{x_\alpha},\frac{x_{\alpha+1}}{x_\alpha},\ldots,\frac{x_n}{x_\alpha}\right)$. We then can write $\CP^{n-1}=(\sqcup^{n}_{\alpha=1}
\varphi_{\alpha}(\U_\alpha))/\simeq$, where, for example,
$$\varphi(\U_1)\ni (z_1,\dots,z_{n-1}) \simeq \left(\frac{1}{z_1},\frac{z_2}{z_1},\ldots, \frac{z_{n-1}}{z_1}\right)\in \varphi_2(\U_2).$$
\indent After formally identifying $\CP^1=\C\cup\{\infty\}$, we define for each positive integer $k$ the $k$-twisted bundle
$$\F^{n}_k=\left(\sqcup^{n}_{\alpha=1}(\U_\alpha\times\CP^1)\right)/\sim$$
\noindent where $\U_\al\times \CP^1\ni ([x_1:\cdots:x_{n}];\xi)\sim ([y_1:\cdots:y_n],\eta)\in U_\be\times\CP^1$ if, and only if, $[x_1:\cdots:x_n]=[y_1:\cdots:y_n]$
and $\eta=\left(\frac{y_\be}{x_\al}\right)^k\xi$. Moreover, let $S_0=\{[x_1:\cdots:x_n];0\}$ and $S_{\infty}=\{[x_1:\cdots:x_n];\infty\}$, which are two global sections
of $\F^{n}_k$. Set $\widehat{\F^{n}_k}=\F^{n}_k\backslash(S_0\cup S_\infty)$ and let $\Psi:\Cnzero \longrightarrow \widehat{\F^{n}_k}$ be the map
$$\Psi:(x_1,x_2,\ldots,x_n)\mapsto ([x_1:\cdots:x_n];x_\al^k)$$
\noindent if $x_\al\neq0$. Because $([x_1:\cdots:x_n];x_\al^k)\sim [x_1:\cdots:x_n];x_\be^k)$ whenever $x_\al\neq0$ and $x_\be\neq0$, the map $\Psi$ is well defined. Furthermore,
one check that $\Psi$ is $k$-to-one and surjective.\\
\indent In particular, when $k=1$, $\Psi$ is one-to-one so $\F^{n}_k=M_n$ can be thought as $\Cnzero$ with $\CP^{n-1}$'s glued at $0$ (the section $S_0$) and $\infty$ (the section $S_\infty$). We next look for suitable metrics on $\Cnzero$ that can be extended to $M_n$.

\subsection{U(2)-Invariant \Kahler\, Metrics on $\Czero$} Let $g$ be a $U(2)$-invariant \Kahler\, metric on $\Czero$, the later with complex coordinates $z=(z^1,z^2)$. Define $u=|z^1|^2$, $v=|z^2|^2$ and $w=u+v$.\\
\indent Since $g$ is a \Kahler\, metric and the de Rham cohomology group $H^2(\Czero)=0$, by the $\partial\bar{\partial}$-lemma one can find a global real smooth function $P:\Czero\longrightarrow\R$ such that\footnote{We are following the notation in \cite[Chapter 2]{Chow}.}
\begin{equation}\label{eq: kahler potential}
g_{\ab}=\dalphabetabar P.\,
\end{equation}
\indent The further assumption of $g$ being rotationally symmetric allows us to write $P=P(r)$, where $r=\log w$. We then set $\varphi(r)=P_r(r)$ (we use subscript for the derivative since later $P$ will be regarded as a function of time as well.) and compute from \ref{eq: kahler potential}
\begin{equation}\label{eq:g}
g=[e^{-r}\varphi\delta_{\alpha\beta} + e^{-2r}(\varphi_r-\varphi)\bar{z}^{\al}z^{\beta}]dz^{\al}d\bar{z}^{\beta}.
\end{equation} 
\noindent Moreover\footnote[2]{The matrix $(g)$ is actually a $4\times 4$ matrix; $(g)=\left(
\begin{matrix}
A & 0\\
0& A\\
\end{matrix}
\right)$, where $A= 
\left(\begin{matrix}
g_{1\bar{1}} & g_{1\bar{2}}\\
g_{2\bar{1}} & g_{2\bar{2}}\\
\end{matrix}
\right)$.},
\begin{equation*}
\left(
\begin{matrix}
g_{1\bar{1}} & g_{1\bar{2}}\\
g_{2\bar{1}} & g_{2\bar{2}}\\
\end{matrix}
\right)=\dfrac{1}{w^2}\left(
\begin{matrix}
v\varphi+u\varphi_r & (\varphi_r-\varphi)\bar{z}^{1}z^{2}\\
(\varphi_r-\varphi)z^{1}\bar{z}^{2} & u\varphi+v\varphi_r\\
\end{matrix}
\right)
\end{equation*}
and thus we find that $\det (g_{\ab})= e^{-2r}\varphi\varphi_r$. From these two past relations we can rapidly observe that a potential $P$ on $\Czero$ gives rise to a \Kahler\, metric as in \ref{eq: kahler potential} if, and only if, 
\begin{equation}\label{eq: Kahler condition}
\varphi>0\, \text{and}\, \varphi_r>0 .    
\end{equation}
\indent For computing the Ricci curvature, recall that on \Kahler\, manifolds the complex Ricci tensor $R_{\ab}dz^\al dz^{\bebar}$ is given by 
\begin{equation*}
R_{\ab}=-\dalphabetabar\log\det g.
\end{equation*}
\indent To put this expression in a more useful form we follow \cite{K} and define
\begin{eqnarray}
G&=&-\log\det g = 2r - \log\varphi-\log\varphi_r\label{eq:principal2}\\
\psi&=&G_r=2-\dfrac{\varphi_r}{\varphi}-\dfrac{\varphi_{rr}}{\varphi_r}\label{eq:principal0}.
\end{eqnarray}
\noindent One next calculates
\begin{equation}\label{eq:Ricci}
R_{\ab}=e^{-r}\psi\delta_{\alpha\beta}+e^{-2r}(\psi_r-\psi)\bar{z}^{\al}z^{\be}
\end{equation}
and thus Ricci can be represented as\footnote[3]{The matrix $(Rc)$ is actually a $4\times 4$ matrix; $(Rc)=\left(
\begin{matrix}
R & 0\\
0& R\\
\end{matrix}
\right)$, where $R= 
\left(\begin{matrix}
R_{1\bar{1}} & R_{1\bar{2}}\\
R_{2\bar{1}} & R_{2\bar{2}}\\
\end{matrix}
\right)$.}
\begin{equation*}
\left(
\begin{matrix}
R_{1\bar{1}} & R_{1\bar{2}}\\
R_{2\bar{1}} & R_{2\bar{2}}\\
\end{matrix}
\right)=\dfrac{1}{w^2}\left(
\begin{matrix}
v\psi+u\psi_r & (\psi_r-\psi)\bar{z}^{1}z^{2}\\
(\psi_r-\psi)z^{1}\bar{z}^{2} & u\psi+v\psi_r\\
\end{matrix}
\right).
\end{equation*}

\indent Also as in \cite{K} we will say that a $(1,0)$ tensor $W$ is an eigenvector of the complex Ricci tensor corresponding to the eigenvalue $\lambda$ if $R_{\ab} W^{\al}=\lambda g_{\ab} W^{\al}.$ Hence understood, the eigenvalues of Rc are
$$
\begin{array}{rcccl}
\lambda_1&=&\dfrac{\psi}{\varphi}&\text{with eigenvector}& U=\bar{z}^2\dfrac{\partial}{\partial z^1}+\bar{z}^1\dfrac{\partial}{\partial z^1}\\
\lambda_2&=&\dfrac{\psi_r}{\varphi_r}&\text{with eigenvector}& V=z^1\dfrac{\partial}{\partial z^2}+z^2\dfrac{\partial}{\partial z^2}.
\end{array}  
$$
\subsection{A Metric of Non-Negative Ricci Curvature} To find a metric of non-negative Ricci curvature on a Calabi line bundle $\CP^1\hookrightarrow M \rightarrow\CP^1$, we wish for a $U(2)$-invariant \Kahler\, potential $P$ on $\Czero$ with the following properties:
\begin{enumerate}
\item $\varphi>0$ everywhere;
\item $\varphi_r>0$ everywhere;
\item $\psi>0$ everywhere;
\item $\psi_r\geq 0$ everywhere;
\item $g$ extends smoothly to a complete metric by adding a $\CP^1$ at $z=0$ ; and
\item $g$ extends smoothly to a complete metric by adding a $\CP^1$ at $z=\infty$.
\end{enumerate} 
\indent To find such $P$, let's rewrite equation \ref{eq:principal0} as 
\begin{equation}\label{principal}
[\log(\varphi\varphi_r)]_r=a,
\end{equation}
\noindent where $a=2-\psi$. We next observe that we can formally solve the ODE above and find formal solution ($I(f)$ stands for an antiderivate of $f$)
\begin{equation}
\varphi(r)=\sqrt{I\left(2e^{I(a)}\right)}. 
\end{equation}
\indent In particular, if $a$ is a constant function, we find that $\varphi(r)=\sqrt{\frac{2}{a}e^{ar}+c}$ is a solution for \ref{principal}. This ($a\equiv 1$) in fact was 
the $\varphi$ used in \cite{K} and, as shown in this particular reference, it leads to potential $P$ satisfying conditions $(1)$ through $(5)$ above, but not $(6)$.\\
\indent We consider instead a function $a$ that is constant equal to 1 in a neighborhood of $r=-\infty$, constant equal to $-1$ in neighborhood of $r=+\infty$, 
and smooth on the whole real line.\\

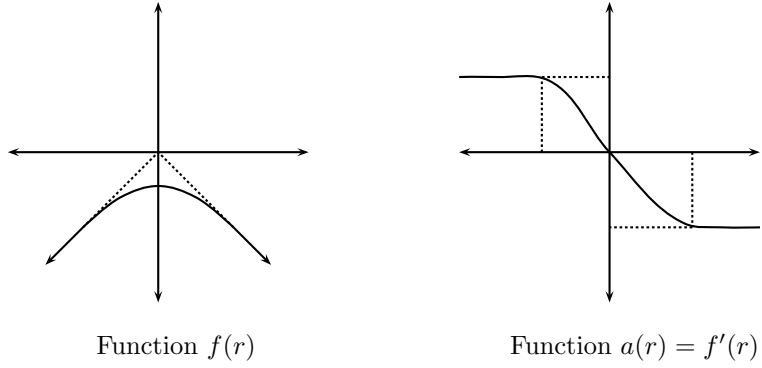
\begin{center}
\psset{unit=1cm}
\begin{pspicture}(-2,-3)(8,2)

\psline{<->}(-2,0)(2,0)\psline{<->}(0,-2)(0,2)
\psline[linestyle=dashed,dash=1pt 1pt](-1,-1)(0,0)
\psline[linestyle=dashed,dash=1pt 1pt](1,-1)(0,0)
\pscurve(-1,-1)(-.5,-.6)(0,-.45)(.5,-.6)(1,-1)
\psline{<-}(-1.5,-1.5)(-1,-1)\psline{<-}(1.5,-1.5)(1,-1)
\uput[0](-1,-2.6){Function $f(r)$}

\psline{<->}(4,0)(8,0)\psline{<->}(6,-2)(6,2)
\pscurve(4,1)(4.5,1)(5.2,.95)(6,0)(7,-.95)(7.5,-1)(8,-1)
\psline[linestyle=dashed,dash=1pt 1pt](5.1,0)(5.1,1)
\psline[linestyle=dashed,dash=1pt 1pt](7.1,0)(7.1,-1)
\psline[linestyle=dashed,dash=1pt 1pt](7.1,-1)(6,-1)
\psline[linestyle=dashed,dash=1pt 1pt](5.1,1)(6,1)
\uput[0](4.5,-2.6){Function $a(r)=f'(r)$}
\end{pspicture}\end{center}
\medskip

\indent To be precise, consider the function $y(r)=-|r|$. Because $y''(r)=\delta_0(r)$, it is possible to smooth $y$ in a small neighborhood of $r=0$ and obtain a smooth function $f(r)$ such that $f''(r)\leq 0$. Then, set $a(r)=f'(r)$; the graph of the function $a$ is shown on the right above.\\
\indent Since 
$$\int_{\R} e^{f(x)}\dx< \int_{\R} e^{-|x|}\dx<+\infty,$$
\noindent the function $F(r)=\int^{r}_{-\infty}e^{f(x)}\dx$ is well defined and smooth. Thus, we find a solution to equation \ref{principal} for $a$ chosen as above:
\begin{equation}\label{solution}
\varphi(r)=\sqrt{2F(r)+c},
\end{equation}
\noindent where $c$ is any arbitrary positive constant (since $F(r)>0$ everywhere). We claim $\varphi$ as above gives a potential $P$ as wished.\\
\indent First, it is clear that $\varphi$ is positive everywhere, and since $\varphi\varphi_r=e^f>0$, $\varphi_r$ is also positive everywhere. Thus $\varphi$ satisfies properties $(1)$ and $(2)$ above, i.e., $P$ gives rise to a \Kahler\, metric.\\
\indent Moreover, $\psi=2-a$, so $\psi>0$ everywhere - since the function $a$ takes values on the interval $[-1,1]$ only. Also (since $f$ concave) $a$ is non-increasing, so $\psi$ is non-decreasing and thus $\psi_r\geq 0$. Hence properties $(3)$ and $(4)$ are also satisfied, that is, $g$ has non-negative Ricci curvature.  \\
\indent To show that $\varphi$ satisfies conditions $(5)$ and $(6)$ it is convenient to compute $F(r)$ when $r$ is near $\pm \infty$ and write explicitly:
\begin{eqnarray}
\varphi=\sqrt{2e^r+c},&&\text{for $r$ near $-\infty$}\label{eq:minusinf};\\
\varphi=\sqrt{-2e^{-r}+c+F_{\infty}},&&\text{for $r$ near $+\infty$\label{eq:plusinf},}
\end{eqnarray}
\noindent where $F_\infty=\int_{\R} e^{f(x)}\dx$ is a constant.\\
\indent Calabi's lemma \cite[Section 3]{C} tells us that $g$ will extend to a smooth \Kahler\, metric on Calabi's twisted line bundle $\CP^1\hookrightarrow M \rightarrow\CP^1$ if $\varphi$ satisfies the following asymptotic properties:
\begin{enumerate}
\item[(i)] There exists positive constants $a_0$ and $a_1$ such that $\varphi$ has the expansion 
$$\varphi(r)=a_0+a_1w+a_2w^2+\mathcal{O}(|w|^3)$$
\noindent as $r\rightarrow -\infty$;
\item[(ii)] There exists a positive constant $b_0$ and a negative constant $b_1$ such that $\varphi$ has the expansion
$$\varphi(r)=b_0+b_1w^{-1}+b_2w^{-2}+\mathcal{O}(|w|^{-3})$$
\noindent as $r\rightarrow \infty$.
\end{enumerate}

\indent For $c>0$, by \ref{eq:minusinf}, $\varphi$ admits the expansion near $r=-\infty$
$$\varphi=\sqrt{c}+\dfrac{1}{\sqrt{c}}w-\dfrac{1}{2c^{3/2}}w^2+
\mathcal{O}(|w|^3).$$
\noindent Thus, by Calabi's lemma, property $(5)$ is verified. And near $r=\infty$, by \ref{eq:plusinf}, $\varphi$ admits the expansion
$$\varphi=\sqrt{c+F_\infty}-\dfrac{1}{\sqrt{c+F_\infty}}w^{-1}-\dfrac{1}
{2(c+F_\infty)^{3/2}}w^{-2}+\mathcal{O}(|w|^{-3}),$$
\noindent so Calabi's lemma, again, guarantees that $g$ satisfy property $(6)$.

\section{Effect of the \Kahler-Ricci Flow}
We finally consider the \Kahler-Ricci Flow evolution $(M,g(t))$, where the initial \Kahler\, potential $P(r,0)$ is taken as the \Kahler\, potential $P(r)$ constructed on the previous section.\\
\indent Because $g$ is rotationally symmetric, let's work at a point $(\zeta,0)\in \Czero$, where $\zeta\neq 0$ is arbitrary. Since
\begin{equation*}
g|_{(\zeta,0)}=\dfrac{1}{|\zeta|^2}\left(\begin{matrix}
\varphi_r& 0\\
0& \varphi
\end{matrix}
\right)
\end{equation*} 
\noindent and 
\begin{equation*}
\text{Rc}|_{(\zeta,0)}=\dfrac{1}{|\zeta|^2}\left(\begin{matrix}
\psi_r& 0\\
0& \psi
\end{matrix}
\right)
\end{equation*} 
\noindent in the standard basis, the \Kahler-Ricci flow equation
$$\dfrac{\partial}{\partial t}g = - \text{Rc}$$
\noindent will be satisfied if, and only if, $\varphi$ evolves by $\varphi_t=-\psi$. This last is reduced to
\begin{equation}\label{eq:evolphi}
\varphi_t=\dfrac{\varphi_{rr}}{\varphi_r}+\dfrac{\varphi_r}{\varphi} -2.
\end{equation}
\indent Moreover, the function $\psi$ (since $\frac{\partial}{\partial r}$ and $\frac{\partial}{\partial r}$ commute) must evolve by
\begin{equation}\label{eq:evolpsi}
\psi_t=\dfrac{\psi_{rr}}{\varphi_r}-\dfrac{\varphi_{rr}\psi_r}{\varphi_r^{2}}
+\dfrac{\psi_r}{\varphi}-\dfrac{\varphi_r\psi}{\varphi^2}.
\end{equation}
\noindent But near $r=-\infty$ we have $\psi(\cdot,0)\equiv 1$, and this reduces equation \ref{eq:evolpsi} to 
\begin{equation}
\dfrac{\partial}{\partial t}\psi\Bigg|_{t=0}=-\dfrac{\varphi_r}{\varphi^2}
\end{equation}
\noindent at the initial time. A final computation then shows that
\begin{equation}
\dfrac{\partial}{\partial t}\psi_r\Bigg|_{t=0}=\dfrac{e^r}{\varphi^5}(e^r-c),
\end{equation} 
\noindent which is strictly negative for $r<\log c$. Hence, the complex Ricci tensor must acquire a negative eigenvalue $\lambda_2<0$ at all points close enough to $r=-\infty$. This completes the proof of Theorem \ref{theorem}. 

\section{Complex Dimension $n\geq 2$}
\indent On this section we expand our result to higher complex dimensions by proving Theorem \ref{theorem2}.\\
\indent If $P$ is now a $U(n)$-invariant \Kahler\, potential on $\Cnzero$, we still have formulas \ref{eq: kahler potential} and \ref{eq:g}. The determinant of the matrix $(g_{\ab})$ is now $\det(g_{\ab})=e^{-nr}\varphi^{n-1}\varphi_r$, and, in particular, conditions in \ref{eq: Kahler condition} are still necessary and sufficient. Formula \ref{eq:principal2} then becomes\footnote[4]{Dan Knopf has drawn to my attention a small oversight on a corresponding formula for $G$ in section 5 of \cite{K}.  As a consequence, few gradient terms are missing on subsequent formulae of \cite{K}, yet these missing terms do not harm the argument made there.}
 
\begin{equation}
G=nr-(n-1)\log\varphi-\log\varphi_r, 
\end{equation}        
so \ref{eq:principal0} is replaced by
\begin{equation}
\psi=G_r=n-(n-1)\dfrac{\varphi_r}{\varphi}-\dfrac{\varphi_{rr}}{\varphi_r},
\end{equation}
\noindent but \ref{eq:Ricci} remains the same. By rotational symmetry, we need to display the eigenvalues of Ricci only at points of the form  $Z=(\zeta,0,0,\ldots,0)\in \Cnzero$, with $\zeta\neq 0$ arbitrary. So we calculate

$$
g|_Z=\dfrac{1}{|\zeta|^2}
\left(\begin{array}{cccc}
\varphi_r &         &       & \\
          & \varphi &       & \\
          &         & \ddots& \\
          &         &       & \varphi          
\end{array}
\right)
$$
\noindent and 
$$
\Rc|_Z=\dfrac{1}{|\zeta|^2}=
\left(\begin{array}{cccc}
\psi_r &         &       & \\
          & \psi &       & \\
          &         & \ddots& \\
          &         &       & \psi          
\end{array}
\right).
$$
\indent We next take $\varphi$ as in \ref{solution}. Hence $\varphi$ and $\varphi_r$ are both greater than zero everywhere. Moreover, we know have $\psi\equiv n-a$, so $\psi>0$ and, since $a$ is non-increasing, $\psi_r\geq 0$,.
\indent Furthermore, expansions of $\varphi$: 
$$\varphi=\sqrt{c}+\dfrac{1}{\sqrt{c}}w-\dfrac{1}{2c^{3/2}}w^2+
\mathcal{O}(|w|^3),$$
\noindent near $r=-\infty$ and 
$$\varphi=\sqrt{c+F_\infty}-\dfrac{1}{\sqrt{c+F_\infty}}w^{-1}-\dfrac{1}
{2(c+F_\infty)^{3/2}}w^{-2}+\mathcal{O}(|w|^{-3}),$$
near $r=\infty$, remain the same so, by Calabi's lemma, $\varphi$ gives raise to a \Kahler\, metric on $\Cnzero$ that can be extended smoothly to $M_n$. \\
\indent We next settle the evolution equations of $\varphi$ and $\psi$:
$$\varphi_t=\dfrac{\varphi_{rr}}{\varphi_r}+(n-1)\dfrac{\varphi_r}{\varphi}-n,
$$
and $$\psi_t=\dfrac{\psi_{rr}}{\varphi_r}-\left(\dfrac{\varphi_{rr}}
{\varphi_r^2}-\dfrac{n-1}{\varphi}\right)\psi_r-(n-1)\dfrac{\varphi_r\psi}
{\varphi^2},$$
\noindent and, since $\psi(\cdot,t)\equiv n-1$ on a neighborhood of $r=-\infty$, this last reduces to:
$$\psi_t=-(n-1)^2\dfrac{\varphi_r}
{\varphi^2}.$$
A further computation then shows that
$$
\dfrac{\partial}{\partial t}\psi_r\Bigg|_{t=0}=(n-1)^2\dfrac{e^r}{\varphi^5}(e^r-c),
$$
\noindent which, as before, is strictly negative for $r<\log c$. This concludes the proof of Theorem \ref{theorem2}.

\begin{remark}
More generally, for each $1\leq k\leq n-1$, one can construct \Kahler\, metrics with non-negative Ricci curvature on $k$-twisted bundles $\F^{n}_{k}$ that immediatly acquire mixed Ricci sign under Ricci Flow. This is achieved by considering the potential $\varphi$ obtained when using the function $ka(r)$ instead of $a(r)$ and then using Calabi's lemma to glue back sections $S_0$ and $S_\infty$.  
\end{remark}

\section{Acknowledgments}
\indent I know of no adequate way of thanking my adviser Dan Knopf for his kind gui\-dan\-ce and unwavering support. I am also thankful for the support received by the National Science Foundation under the agreement No. DMS-0545984. Finally, I extend my gratitude to Carlo Mantegazza and Wolfgang Ziller for support received during visits to the Centro di Ricerca Matematica Ennio De Giorgi and to the Ins\-ti\-tu\-to Nacional de Matem\'atica Pura e Aplicada (IMPA), respectively, and to Huai-Dong Cao for helpful comments on the original manuscript.\\
\bibliographystyle{amsplain}

\begin{thebibliography}{99}

\bibitem{B}
{\bf Bando, Shigetoshi.} 
\newblock On the classification of three-dimensional compact \Kahler\, manifolds of non-negative bisectional curvature. 
\newblock {\it J. Differential Geom.} {\bf 19}(2) (1984) 283--297.

\bibitem{BW}
{\bf B\"{o}hm, Christoph; Wilking, Burkhard.} 
\newblock Nonnegatively curved manifolds with finite fundamental groups admit metrics with positive Ricci curvature.  
\newblock {\it Geom. Funct. Anal.}  {\bf 17}  (2007),  no. 3, 665--681.

\bibitem{BW1}
{\bf B\"{o}hm, Christoph; Wilking, Burkhard.} 
\newblock Manifolds with positive curvature operators are space forms. \newblock {\it Ann. of Math. (2)} {\bf  167}  (2008),  no. 3, 1079--1097.

\bibitem{Bre}
{\bf Brendle, Simon.} 
\newblock A general convergence result for the Ricci flow in higher dimensions.  
\newblock {\it Duke Math. J.}  {\bf 145}  (2008),  no. 3, 585--601.

\bibitem{BS}
{\bf Brendle, Simon; Schoen, Richard.} 
\newblock {Manifolds with $1/4$-pinched curvature are space forms.} 
\newblock {\it J. Amer. Math. Soc.} {\bf 22}  (2009),  no. 1, 287--307.

\bibitem{C}
{\bf Calabi, Eugenio.}
\newblock Extremal \Kahler$\,$metrics,
\newblock {\it Seminar on Differential Geometry},  pp. 259--290, Ann. of Math. Stud., {\bf 102}, Princeton Univ. Press, Princeton, N.J., 1982.
53C55 (58E30) 

\bibitem{Cao1}
{\bf Cao, Huai-Dong.} 
\newblock Existence of gradient \Kahler-Ricci solitons.
\newblock {\it Elliptic and parabolic methods in geometry} (Minneapolis, MN, 1994), 1--16, A K Peters, Wellesley, MA, 1996. 

\bibitem{Cao2}
{\bf Cao, Huai-Dong.}  
\newblock Limits of solutions to the \Kahler-Ricci flow.
\newblock {\it J. Differential Geom.} {\bf 45} (1997), no. 2, 257--272. 

\bibitem{Cao3}
{\bf Cao, Huai-Dong.}
\newblock Deformation of \Kahler\, metrics to \Kahler-Einstein metrics on compact \Kahler\, manifolds.
\newblock {\it Invent. Math.} {\bf 81} (1985), no. 2, 359�372.

\bibitem{CaoH}
{\bf Cao, Huai-Dong; Hamilton, Richard S.} 
\newblock On the preservation of positive orthogonal bisectional curvature under the \Kahler-Ricci flow.
\newblock {\it private communication} (unpublished work) (1988)

\bibitem{HC}
{\bf Chen, Haiwen.}
\newblock {Pointwise $\frac14$-pinched $4$-manifolds.}  
\newblock {\it Ann. Global Anal. Geom.}  {\bf  9}  (1991),  no. 2, 161--176.

\bibitem{XC}
{\bf Chen, Xiuxiong.} 
\newblock On \Kahler\, manifolds with positive orthogonal bisectional curvature.  
\newblock {\it Adv. Math.}  {\bf 215}  (2007),  no. 2, 427--445. 

\bibitem{XCL}
{\bf Chen, Xiuxiong; Li, Haozhao.} 
\newblock The \Kahler-Ricci flow on \Kahler\, manifolds with 2-non-negative traceless bisectional curvature operator. 
\newblock {\it Chin. Ann. Math. Ser.} B {\bf 29}  (2008),  no. 5, 543--556.

\bibitem{XCT1}
{\bf Chen, Xiuxiong; Tian, Gang.}
\newblock Ricci flow on \Kahler-Einstein surfaces. 
\newblock {\it Invent. Math.} {\bf 147} (2002), no. 3, 487--544. 

\bibitem{XCT2}
{\bf Chen, Xiuxiong; Tian, Gang.}
\newblock Ricci flow on \Kahler-Einstein manifolds.  
\newblock {\it Duke Math. J.}  {\bf 131}  (2006),  no. 1, 17--73.


\bibitem{CNL}
{\bf Chow, Bennett; Lu, Peng; Ni, Lei.} 
\newblock {\it Hamilton's Ricci flow.} 
\newblock Graduate Studies in Mathematics, {\bf 77}. American Mathematical Society, Providence, RI; Science Press, New York, 2006. 

\bibitem{Chow}
{\bf Chow, Bennett; Chu, Sun-Chin; Glickenstein, David; Guenther, Christine; Isenberg, James; Ivey, Tom; Knopf, Dan; Lu, Peng; Luo, Feng; Ni, Lei}
\newblock {\it The Ricci Flow: Techniques and Applications, Part I: Geometric Aspects.} 
\newblock Mathematical Surveys and Monographs, {\bf 135}. American Mathematical Society, Providence, RI, 2007.

\bibitem{FIK}
{\bf Feldman, Mikhail; Ilmanen, Tom; Knopf, Dan.}
\newblock Rotationally symmetric shrinking and expanding gradient \Kahler-Ricci solitons.
\newblock {\it J. Differential Geom.} {\bf 65} (2003), no. 2, 169--209. 

\bibitem{GK}
{\bf Goldberg, Samuel I.; Kobayashi, Shoshichi.}
\newblock Holomorphic bisectional curvature.
\newblock {\it J. Differential Geometry} {\bf 1} 1967 225--233. 

\bibitem{H1}
{\bf Hamilton, Richard S.} 
\newblock {Three-manifolds with positive Ricci curvature}. 
\newblock {\it J. Differential Geom.} {\bf 17} (1982), no. 2, 255--306.

\bibitem {H2}
{\bf Hamilton, Richard S.} 
\newblock Four-manifolds with positive curvature operator.  
\newblock {\it J. Differential Geom.}  {\bf 24}  (1986),  no. 2, 153--179.

\bibitem{H3}
{\bf Hamilton, Richard S.} 
\newblock Four-manifolds with positive isotropic curvature.  
\newblock {\it Comm. Anal. Geom.}  {\bf 5}  (1997),  no. 1, 1--92.

\bibitem{K}
{\bf Knopf, Dan.}
\newblock Positivity of Ricci curvature under the \Kahler-Ricci flow.
\newblock {\it Commun. Contemp. Math.} {\bf 8} (2006), no. 1, 123--133. 

\bibitem{Ngu}
{\bf Nguyen, Huy T.}
\newblock Isotropic Curvature and the Ricci Flow.
\newblock {\it Int Math Res Notices}. (2009) to appear.


\bibitem{Ni}
{\bf Ni, Lei.}
\newblock Ricci flow and nonnegativity of sectional curvature. 
\newblock {\it Math. Res. Lett.}  {\bf 11}  (2004),  no. 5-6, 883--904.

\bibitem{Mok}
{\bf Mok, Ngaiming.}
\newblock The uniformization theorem for compact \Kahler\, manifolds of non-negative holomorphic bisectional curvature.
\newblock {\it J. Differential Geom.} {\bf 27}(2) (1988) 179--214.

\bibitem{PS1}
{\bf Phong, D. H.; Sturm, Jacob.} 
\newblock On the \Kahler-Ricci flow on complex surfaces.
\newblock  {\it Pure Appl. Math. Q.}  {\bf 1}  (2005),  no. 2, part 1, 405--413.

\bibitem{PS2}
{\bf Phong, D. H.; Sturm, Jacob.}
\newblock On stability and the convergence of the \Kahler-Ricci flow. 
\newblock {\it J. Differential Geom.} {\bf 72} (2006), no. 1

\bibitem{Sh1} 
{\bf Shi, Wan-Xiong.} 
\newblock Complete noncompact three-manifolds with nonnegative Ricci curvature. 
\newblock {\it J. Differential Geom.} {\bf 29} (1989), 353--360. 

\bibitem{Sh2}
{\bf Shi, Wan-Xiong.} 
\newblock Deforming the metric on complete Riemannian manifolds. 
\newblock {\it J. Differential Geom.} {\bf 30} (1989), 223--301. 

\bibitem{Sh3}
{\bf Shi, Wan-Xiong.} 
\newblock Ricci flow and the uniformization on complete noncompact \Kahler\, manifolds. 
\newblock {\it J. Differential Geom.} {\bf 45} (1997), 94--220.  

\bibitem{W}
{\bf Wilking, Burkhard.}
\newblock Sharp estimates for the Ricci flow.
\newblock{\it International Symposium on Differential Geometry ``In honor of Marcos Dajczer on his 60th birthday"} 2009.
\end{thebibliography}

\end{document}